\documentclass[12pt]{amsart}
\usepackage{graphicx}
\usepackage{a4wide}
\usepackage[centertags]{amsmath}
\usepackage{amsfonts}
\usepackage{amssymb}
\usepackage{amsthm}
\usepackage{hyperref}
\newtheorem{thm}{Theorem}%[section]
\newtheorem{lem}{Lemma}[section]
\newtheorem{cor}{Corollary}

\newtheorem{prop}[lem]{Proposition}
\theoremstyle{definition}

\theoremstyle{remark}
\newtheorem{rem}{Remark}[section]
\numberwithin{equation}{section}

\newcommand{\set}[1]{\left\{#1\right\}}

\newcommand{\calG}{\mathcal{G}}

\newcommand{\calK}{\mathcal{K}}

\newcommand{\bbZ}{\mathbb{Z}}

\newcommand{\bbR}{\mathbb R}

\newcommand{\bbN}{\mathbb N}

\newcommand{\SU}{ \mathrm{SU}}

\begin{document}
\title[Distribution modulo prime powers]
{Distribution of twisted Kloosterman sums\\ modulo prime powers}%
\author{Dubi Kelmer }%
\address{School of Mathematics,
Institute for Advanced Study, 1 Einstein Drive , Princeton, New Jersey 08540 US}
\email{kelmerdu@ias.edu}

\thanks{
This material is based upon work supported by the National Science Foundation under agreement No. DMS-0635607. Any opinions, findings and conclusions or recommendations expressed in this material are those of the author and do not necessarily reflect the views of the National Science Foundation.}%
%\subjclass{11L05}%
%\keywords{Kloosterman sums, equidistribution}%

\date{\today}%
\dedicatory{}%
\commby{}%

\begin{abstract}
In this note we study Kloosterman sums twisted by a multiplicative characters
modulo a prime power. We show, by an elementary calculation, that these sums become equidistributed on the real
line with respect to a suitable measure.
\end{abstract}
\maketitle
\section{Introduction}
Let $p$ be a prime and $q=p^k$ a prime power.
For any multiplicative character $\chi\in\widehat{(\bbZ/q\bbZ)^*}$, and any $a,b\in(\bbZ/q\bbZ)^*$ consider the twisted Kloosterman sum
\[K_q(a,b,\chi)=\sum_{x\in(\bbZ/q\bbZ)^*}e_q(ax+bx^{-1})\chi(x),\]
where $e_q(x)=\exp(\frac{2\pi i x}{q})$ is an additive character of $\bbZ/q\bbZ$.

When $\chi=1$ is trivial, the sum $K_q(a,b,1)=K_q(ab,1,1)\in\bbR$ is real and satisfies $|K_q(a,b,1)|\leq 2\sqrt{q}$.
For $q=p$ prime this is the Weil bound \cite{Weil48}, and for prime powers $q=p^k$ this can be shown by an elementary calculation (see Corollary \ref{c:bound}).
In \cite{Katz88}, following Deligne's Equidistribution Theorem for Frobenius conjugacy classes, Katz showed that as $p\to\infty$ through primes, the normalized sums $p^{-\frac{1}{2}}K_p(a,1,1)$ (with $a$ varying over $(\bbZ/p\bbZ)^*$) become equidistributed on the real line with respect to the Sato-Tate measure $$\mu_{\mathrm{ST}}(f)=\frac{1}{2\pi}\int_{-2}^2f(x)\sqrt{4-x^2}dx.$$

One can also consider the distribution of the twisted sums when the varying parameter is the multiplicative character.
To be more precise, we consider the family $\calK_q(a,\chi)=q^{-\frac{1}{2}}K_q(a,-a,\chi)$ with $a\in\bbZ$ fixed and $\chi$ varying over $(\bbZ/q\bbZ)^*$ (we take $b=-a$ to insure that the sums are real and the normalization is chosen to get mean square one).
This sort of question is similar to the distribution of exponential sums considered by Kurlberg and Rudnick with respect to the matrix elements of the quantum cat map \cite{KurRud05}.

For $q=p$ a prime, we expect the sums $\calK_p(a,\chi)$ to become equidistributed (as $p\to\infty$) with respect to the Sato-Tate measure.
Note that since the multiplicative character is not an algebraic parameter, such an equidistribution result does not follow from Deligne's Equidistribution Theorem as in the case of the sums $K_p(a,1,1)$. Nevertheless, this result should follow from Katz's analysis of Tannakian categories \cite{Katz03}. (We also remark that the numerical evidence for this is quite convincing.)

In this note we study the distribution of these normalized sums for $q=p^k$ a (nontrivial) prime power.
We show that as $p\to\infty$, for any fixed $a\in\bbZ$, the normalized sums $\calK_q(a,\chi)$ also become equidistributed on the real line but with respect to the measure
$$\mu(f)=\frac{1}{2}f(0)+\frac{1}{2\pi}\int_{-2}^2\frac{f(x)dx}{\sqrt{4-x^2}}.$$
Moreover, we show that for any $r$ distinct integers $a_1,\ldots,a_r$, the corresponding $r$-tuples
$(\calK_q(a_1,\chi),\ldots,\calK_q(a_r,\chi))$ become equidistributed in $\bbR^{r}$ with respect to the product measure $d\mu^r$.

\begin{thm}\label{t:main}
Fix $r$ distinct integers $a_1,\ldots,a_r$. For any continuous bounded function $g\in C(\bbR^r)$ and any $k>1$,
\[\lim_{p\to\infty}\frac{1}{p^k}\sum_\chi g(\calK_{p^k}(a_1,\chi),\ldots,\calK_{p^k}(a_r,\chi))=\int_{\bbR^r} g(x) d\mu^r(x).\]
\end{thm}

\begin{rem}
One can also study the limiting distribution of the sums $K_q(a,1,1)$ for $q$ a prime power and $a$ varying over $(\bbZ/q\bbZ)^*$.
For this case the normalized sum is given (up to a sign) by $q^{-1/2}K_q(a,1,1)=\pm 2\cos(\frac{4\pi c}{q})$ if $a\equiv c^2\pmod{q}$ and it vanishes otherwise. From this, one can easily deduce that as $q\to\infty$ these sums also become equidistributed with respect to the measure $\mu$ above. We also note that for Salli\'{e} sums (i.e., for $K_q(a,1,\chi_2)$ with $\chi_2$ the quadratic character) the same formula holds also for prime modulus (as well as for prime powers). Hence, the sums $K_p(a,1,\chi_2)$ (when $a$ varies and $p\to\infty$) also become equidistributed with respect to $\mu$.
\end{rem}

\begin{rem}
There is a natural group theoretic interpretation of the measures $\mu$ and $\mu_{ST}$ defined above. For any subgroup $K\subseteq \SU(2)$, the trace map sends the Haar measure of $K$ to a corresponding measure on $[-2,2]$. The Sato-Tate measure $\mu_{ST}$ is obtained in this manner if we take $K$ to be the full group. The measure $\mu$, on the other hand, is obtained if we take $K$ to be the normalizer of a maximal torus in $\SU(2)$. This interpretation is natural when considering Kloosterman or Salli\'{e} sums with prime modulus, since in this case the corresponding sums can be represented as traces of matrices from $\SU(2)$. Deligne's Equidistribution Theorem them implies that these matrices become equidistributed in the corresponding group (i.e., in $\SU(2)$ for $K(a,1,1)$ and in the normalizer of a torus for $K(a,1,\chi_2)$). For the prime power case we are not aware of such an interpretation of these sums, however, in this case the proof of the limiting distribution is elementary.
\end{rem}

\section*{Acknowledgments}
I thank Nick Katz for discussions on equidistribution of exponential sums. I also thank Emmanuel Kowalski and Ze{\'e}v Rudnick for their comments.

\section{Elementary Calculation}
For exponential sums modulo a prime power $q=p^k$ it is possible to evaluate the sum by elementary means (see e.g., \cite[Section 12.3]{IwaniecKowalski} or \cite[Chapter 1.6]{BerndtEvansWilliams}).
Applying this method specifically for Kloosterman sums gives the following results.
\begin{prop}
Let $q=p^{2l}$ be an even power. Then for any character $\chi\in\widehat{(\bbZ/q\bbZ)^*}$ there is a unique $t_\chi\in\bbZ/p^l\bbZ$ such that
$\chi(1+p^lx)=e_{p^l}(t_\chi x)$. In this case the Kloosterman sum is given by
\[K_q(a,b,\chi)=p^l\mathop{\sum_{x\in(\bbZ/p^l\bbZ)^*}}_{h(x)\equiv 0\pmod{p^l}}e_q(ax+bx^{-1})\chi(x),\]
with $h(x)=ax^2+t_\chi x-b$.
\end{prop}
\begin{proof}
See \cite[lemma 12.2]{IwaniecKowalski}).
\end{proof}
\begin{prop}
Let $q=p^{2l+1}$ be an odd power. Then for any character $\chi\in\widehat{(\bbZ/q\bbZ)^*}$ there is a unique $t_\chi\in\bbZ/p^{l+1}\bbZ$ such that
$$\chi(1+p^lx+p^{2l}\frac{x^2}{2})=e_{p^{l+1}}(t_\chi x).$$ In this case the Kloosterman sum is given by
\[\calK_q(a,b,\chi)=p^l\mathop{\sum_{x\in(\bbZ/p^l\bbZ)^*}}_{h(x)\equiv 0\pmod{p^l}}e_q(ax+bx^{-1})\chi(x)\calG(x),\]
with $h(x)=ax^2+t_\chi x-b$, and $\calG(x)$ the Gauss sum
\[\calG(x)=\sum_{y\pmod{p}}e_p(d(x)y^2+p^{-l}h(x)y)\]
with $d(x)=\frac{(p-1)}{2}t_\chi x^2+bx$

\end{prop}
\begin{proof}
See \cite[lemma 12.3]{IwaniecKowalski}).
\end{proof}

\begin{cor}\label{c:bound}
For any character $\chi$ with $t_\chi^2\not\equiv -4ab \pmod p$ we have $|K_q(a,b,\chi)|\leq 2\sqrt{q}$.
In particular, for any $a,b\in(\bbZ/q\bbZ)^*$ we have $|K(a,b,1)|\leq 2\sqrt{q}$ and for any $a\in(\bbZ/q\bbZ)^*$ and $\chi$ with $t_\chi\not\equiv\pm 2a\pmod{p}$ we have $|K(a,-a,\chi)|\leq 2\sqrt{q}$.
\end{cor}
\begin{proof}
The condition $t_\chi^2\not\equiv -4ab \pmod p$ implies that the polynomial $h(x)=ax^2+t_\chi x+b$ is separable, and hence the equation $h(x)\equiv 0\pmod{p^l}$ has two or zero solutions in $\bbZ/p^l\bbZ$.
When $k=2l$ is even, $q^{-1/2}K(a,b,\chi)=p^{-l}K_q(a,b,\chi)$ is given by a sum of two (or zero) elements each of absolute value one implying the bound.
For $k=2l+1$ odd, we also use the identity $|\calG(x)|=\sqrt{p}$. For this to hold we have to exclude the possibility that $d(x)\equiv 0\pmod{p}$ and $h(x)\equiv 0\pmod{p^{l+1}}$ (in which case $\calG(x)=p$). However the equation
$d(x)=h(x)= 0\pmod{p}$ has a solution only if $t_\chi^2=-4ab\pmod{p}$.
\end{proof}

%*******************************************************************************************************************************************
%********************************* SECTION 3 Equidistribution ******************************************************************************
%*******************************************************************************************************************************************
\section{Equidistribution}
We now turn to the  proof of Theorem \ref{t:main}.
Fix a finite set of $r$ nonzero distinct integers $a_1,\ldots,a_r$.
We need to show that for any continuous bounded function $g\in C(\bbR^r)$
\begin{eqnarray*}
\lefteqn{\lim_{p\to\infty}\frac{1}{q}\sum_{\chi}g(\calK_q(a_1,\chi),\ldots,\calK_q(a_r,\chi))=}\\
&&\int_{\bbR^r}g(x_1,\ldots,x_r)d\mu(x_1)\cdots d\mu(x_r).
\end{eqnarray*}
We will use the fact that the normalized sums are (almost always) bounded by $2$ to reduce this to a moment calculation.
We then exploit the explicit calculation of the sums to reduce the moment calculation to a counting argument.
\subsection{Reduction to a moment calculation}
For every $0\leq j\leq k$ let
\[C_p(k,j)=\set{\chi\in \widehat{(\bbZ/q\bbZ)^*}|\chi(x)=1,\forall x\equiv 1\pmod{p^j}}.\]
For $\chi,\chi'$ any two characters, and $j\leq [\frac{k}{2}]$ (the integer part of $\frac{k}{2}$), we have that $t_{\chi}\equiv t_{\chi'}\pmod{p^j}$ if and only if $\chi'\chi^{-1}\in C_p(k,k-j)$.

For any $q=p^k$ we define the set of characters
\[S_q=S_q(a_1,\ldots,a_r)=\set{\chi\in\widehat{(\bbZ/q\bbZ)^*}|t_\chi\not\equiv \pm 2a_j\pmod{p}},\]
where $t_\chi$ is determined by $\chi$ as above.
Since the order of $C_p(k,j)$ is $p^{j-1}(p-1)$, we get that
for any $a_j$ there are precisely $2p^{k-2}(p- 1)$ characters with $t_\chi\equiv\pm 2 a_j\pmod{p}$. Hence $\frac{|S_q|}{q}=1+O(\frac{1}{p})$ and it is sufficient to show that for any continues bounded function $g\in C_0(\bbR^r)$
\begin{eqnarray*}
\lefteqn{\lim_{p\to\infty}\frac{1}{q}\sum_{\chi\in S_q}g(\calK_q(a_1,\chi),\ldots,\calK_q(a_r,\chi))=}\\
&&\int_{[0,\pi]^d}g(x_1,\ldots,x_r)d\mu(x_1)\cdots d\mu(x_r).
\end{eqnarray*}

Now note that for any $\chi\in S_q$ we have the bound $|\calK_q(a_j,\chi)|\leq 2$. Since the measure $\mu$ is also supported on $[-2,2]$ it is sufficient to check this for continues function on $[-2,2]^r$.
Finally, note that any continues function on $[-2,2]^r$ can be uniformly approximated by polynomials, so it is sufficient to check all mixed moments.
We have thus reduced Theorem \ref{t:main} to the following proposition
\begin{prop}\label{p:moments}
For any  $m_1,\ldots m_r\in\bbN$,
\[
\lim_{p\to\infty}\frac{1}{q}\sum_{\chi\in S_q}\prod_{j}\calK_q(a_j,\chi)^{m_j}= \prod_j\left(\int_{-2}^2 x^{m_j}d\mu(x)\right)
\]
\end{prop}
\subsection{Reduction to a counting argument}
We now reduce Proposition  \ref{p:moments} to a simple counting argument.
We need to show that for any  $m_1,\ldots, m_r\in\bbN$,
\[
\lim_{p\to\infty}\frac{1}{q}\sum_{\chi\in S_q}\prod_{j}\calK_q(a_j,\chi)^{m_j}= \prod_j\left(\int_{-2}^2 x^{m_j}d\mu(x)\right)
\]
With out loss of generality, we can assume all the $m_j$ are nonzero in which case the right hand side is given by
\begin{eqnarray*}
 \prod_j\left(\int_{-2}^2 x^{m_j}d\mu(x)\right)=\prod_j\left(\frac{1}{2}\int_0^\pi \left(2\cos(\theta)\right)^{m_j}\frac{d\theta}{\pi}\right)
\end{eqnarray*}
where we made the change of variables $x=2\cos(\theta)$.
The integrals are easily computed to give
\[\int_0^\pi \left(2\cos(\theta)\right)^{m_j}\frac{d\theta}{\pi}=\left\lbrace\begin{array}{cc}
\begin{pmatrix}m_j\\n_j\end{pmatrix} & m_j=2n_j \mbox{ is even}\\
0 & \mbox{ otherwise}\end{array}\right.\]

Now for the left hand side. Let $l=[\frac{k}{2}]$ (the integer part of $\frac{k}{2}$). For any character $\chi\in S_{q}$ let $t_\chi\in\bbZ/p^{k-l}\bbZ$ as above. Note that if we multiply $\chi$ by any character in $C_p(k,l)$ this does not change $t_\chi$. Consequently, if for any $t\in\bbZ/p^{k-l}\bbZ$ we fix a representative $\chi_t$ (s.t., $t_{\chi_t}=t$), we can write
\begin{eqnarray*}\lefteqn{\frac{1}{q}\sum_{\chi\in S_q}\prod_{j} \calK_q(a_j,\chi)^{m_j}=}\\
&&=\frac{1}{q}\sum_{t\neq \pm a_j}\sum_{\chi\in C_p(k,l)}
\prod_{j}\calK_q(a_j,\chi\chi_t)^{m_j}
\end{eqnarray*}
(where the first sum is only over $t\in\bbZ/p^{k-l}\bbZ$ such that $\forall j,\;t\neq\pm 2a_j\pmod{p}$).

Next note that the only contributions to this sum,
comes from $t\in\bbZ/p^{k-l}\bbZ$ such that for any $j$ there is $x_j\in(\bbZ/p^l\bbZ)^*$ satisfying $a_j(x_j+x_j^{-1})\equiv t\pmod{p^l}$.
For such $t$ it is given by
\[\calK_{q}(a_j,\chi\chi_t)=2\Re(e_{q}(a_j(x_j-x_j^{-1}))\chi\chi_t(x_j) e(\alpha(x_j)))=2\cos(\theta(a_j,\chi\chi_t)),\]
with $x_j\in(\bbZ/p^l\bbZ)^*$ satisfying $a_j(x_j+x_j^{-1})\equiv t\pmod{p^l}$ and $\alpha(x_j)$ is the angle of the normalized Gauss sum $\calG(x_j)$ for $k$ odd and zero for $k$ even.

Let
\[Y(p^l)=\set{\vec{x}\in({(\bbZ/p^l\bbZ)^*})^r|a_1(x_1+x_1^{-1})=a_j(x_j+x_j^{-1}),\;\forall 2\leq j\leq r},\]
and for every $\vec{x}\in Y(p^l)$ let $t(\vec{x})\equiv a_1(x_1+x_1^{-1})\in(\bbZ/p^l\bbZ)^*$. We can thus rewrite the above sum as
\begin{eqnarray*}\lefteqn{\frac{1}{q}\sum_{\chi\in S_q}\prod_{j} \calK_q(a_j,\chi)^{m_j}=}\\
&&=\frac{1}{2^rq}\sum_{Y'(p^l)}\sum_{t\equiv t(\vec{x})}\sum_{\chi\in C_p(k,l)}
\prod_{j}(2\cos(\theta(a_j,\chi\chi_t)))^{m_j}
\end{eqnarray*}
where the middle sum is over the co-set $\set{t\in\bbZ/p^{k-l}\bbZ|t\equiv t(\vec{x})\pmod{p^l}}$
and the notation $Y'(p^l)$ means that we exclude elements with $x_j=\pm 1\pmod{p}$. Notice that for $k$ even the middle sum is just one element and for $k$ odd it is a sum over $p$ elements.

Now use the formula,
\[(2\cos(\theta(a_j,\chi\chi_t)))^{m_j}=\sum_{n=0}^{m_j}\begin{pmatrix} {m_j}\\n\end{pmatrix}\cos((m_j-2n)\theta(a_j,\chi\chi_t)).\]
The main contribution comes from the terms where in all the factors $m_j-2n_j=0$. This vanishes unless all $m_j=2n_j$ are even in which case it is given by
\[\frac{|Y'(p^l)||C_p(k,l)|}{p^{2l}}\prod_{j}\frac{1}{2}\begin{pmatrix} m_j\\ n_j\end{pmatrix}.\]
Since $C_p(k,l)=p^{l}+p^{l-1}$ to get the correct main term we need to show that $|Y'(p^l)|=p^l+o(p^l)$ as $p\to\infty$
 (which is shown in Lemma \ref{l:count1}).

We also have an error term coming from the sums of the form
\[\frac{1}{q}\sum_{\vec{x}\in Y'(p^l)}\sum_{t\equiv t(\vec{x})}\sum_{\chi\in C_p(k,l)}
\prod_{j}\cos(n_j\theta(a_j,\chi\chi_t)),\]
with $\{n_1,\ldots, n_r\}$ nonzero integers.
For any $\vec{x}\in Y'(p^l)$ and $t\equiv a_1(x_1+x_1^{-1})\pmod{p^l}$ we have that
\[\cos(n_j\theta(a_j,\chi\chi_t)=2\Re[e_{q}(n_ja_j(x_j-x_j^{-1}))\chi\chi_t(x_j^{n_j})e(n_j\alpha(x_j))].\]
 Hence, to bound the error term we need to bound sums of the form
\[\frac{1}{q}\sum_{\vec{x}\in Y(p^l)}\mathop{\sum}_{t\equiv t(\vec{x})}\sum_{\chi\in C_p(k,l)}e_{p^k}(\sum_{j=1}^r n_j(a_j(x_j-x_j^{-1})+\alpha(x_j)))\chi\chi_t(\prod_{j=1}^rx_j^{n_j}).\]
Rewrite this sum as
\[\frac{1}{q}\sum_{\vec{x}\in Y(p^l)}e_{p^k}(\sum_{j=1}^r n_j(a_j(x_j-x_j^{-1})+\alpha(x_j)))\sum_{t\equiv t(\vec{x})}\chi_t(\prod_{j=1}^rx_j^{n_j})\sum_{\chi\in C_p(k,l)}\chi(\prod_{j=1}^r x_j^{n_j}),\]
and note that unless  $\prod_{j=1}^r x_j^{n_j}\equiv 1\pmod {p^l}$ the sum $\sum_{\chi\in C_p(k,l)}\chi(\prod_{j=1}^r x_j^{n_j})=0$.
We can thus bound the above sum by $\frac{\sharp C_p(k,l)}
{p^{2l}}=p^{-l}+p^{-(l+1)}$ times the number of elements in
\[Y_0(p^l)=\set{\vec{x}\in Y(p^l)|x_1^{n_1}\cdots x_r^{n_r}\equiv 1\pmod{p^l}}.\]

Now the proof of Proposition \ref{p:moments} is concluded by the following two counting lemmas:
\begin{lem}\label{l:count1}
As $p\to\infty$, the number of points in $Y'(p^l)$ satisfy
\[\sharp Y'(p^l)=p^l+O(p^{l-\frac{1}{2}})\]
\end{lem}
and
\begin{lem}\label{l:count2}
As $p\to\infty$, the number of points in $Y'_0(p^l)$ satisfy
$$\sharp Y'_0(p^l)=O(p^{l-1}).$$
\end{lem}

%***************************************************************************************************************
%***************************************************************************************************************
\subsection{Counting arguments}
We conclude this section with the proof of the two counting lemmas
\begin{proof}[\textbf{Proof of Lemma \ref{l:count1}}]
For any $t\in(\bbZ/p^l\bbZ)^*$ satisfying $\forall j,\;t\neq\pm 2a_j\pmod{p}$ we have that
\[\sharp\set{\vec{x}|\forall j,\;a_j(x_j+x_j^{-1})=t}=\left\lbrace\begin{array}{cc}
2^r & \forall j,\;t^2-4a_j^2\equiv \square\pmod{p}\\
 0 & \mbox{ otherwise}
\end{array}\right.
\]
% If, on the other hand if $t\equiv \pm 2a\pmod{p}$ for some $j$, then
% \[\sharp\set{\vec{x}|\forall j,\;a_j(x_j+x_j^{-1})=t}\leq 2^{r}p^{k-1}\]
% as there are at most two possibilities for $x_i$ with $i\neq j$ and at most $2p^{k-1}$ possibilities for $x_j$.
We thus have
\begin{eqnarray*}\sharp Y'(p^l) &= &\mathop{\sum_{t\in(\bbZ/p^l\bbZ)^*}}_{t\not\equiv\pm 2a_j(p)}\sharp\set{\vec{x}\in((\bbZ/p^l\bbZ)^*)^r|\forall j,\;a_j(x_j+x_j^{-1})=t}\\
&= &2^rp^{k-1}\sharp\set{t\in(\bbZ/p\bbZ)^*| \forall j,\;0\neq t^2-4a_j^2 \equiv \square\pmod{p}}.
\end{eqnarray*}
To conclude the proof we use the estimate
\[2^r\sharp\set{t\in(\bbZ/p\bbZ)^*| \forall j,\;t^2-4a^2 \equiv \square\pmod{p}}=p+O(\sqrt{p}).\]
(To get this estimate write
\[2^r\sharp\set{t\in\bbZ/p\bbZ| \forall j,\;t^2-4a^2 \equiv \square\pmod{p}}=\sum_{t}\prod_{j=1}^r\left(\chi_2(t^2-4a^2)+1\right)\]
with $\chi_2$ the quadratic character modulo $p$, and apply the Weil bounds on the corresponding character sums.)
\end{proof}

\begin{proof}[\textbf{Proof of Lemma \ref{l:count2}}]
To prove the bound  $Y_0(p^l)=O(p^{l-1})$ we will show that there is a nonzero polynomial $F(t)$ with coefficients in $\bbZ[\frac{1}{a_1},\ldots,\frac{1}{a_r}]$,
such that for any $\vec{x}\in Y_0(p^l)$, with $t\equiv a_1(x_1+x_1^{-1})\pmod{p^l}$ we have $F(t)\equiv 0\pmod{p^l}$. This would imply that $t$ can take at most $\deg F$ values modulo $p$, implying that $\sharp Y_0(p^l)\leq 2^r\deg(F)p^{l-1}$.

Now to define $F$, consider the formal polynomial in the variables $x_1^{\pm 1},\ldots x_r^{\pm 1}$ given by
\[G(x_1,\ldots x_r)=\prod_{\sigma\in\{\pm 1\}^r}\left(\prod_{j=1}^r x_j^{\sigma_jn_j}-1\right).\]
This polynomial is symmetric under any substitution $x_j\mapsto x_j^{-1}$ and hence there is another polynomial $\tilde{F}$ in $r$ variables with integer coefficients (of degree bounded by $2^r\max\{|n_j|\}$) satisfying
 \[G(x_1,\ldots,x_r)=\tilde{F}(x_1+x_1^{-1},\ldots,x_r+x_r^{-1}).\]
Define the polynomial $F(t)=\tilde{F}(\frac{t}{a_1},\ldots,\frac{t}{a_r})$. For any $x_1,\ldots,x_r$ with $a_j(x_j+x_j^{-1})=t\pmod{p^l}$ we have
\[G(x_1,\ldots,x_r)=\tilde{F}(\frac{t}{a_1},\ldots,\frac{t}{a_r})=F(t).\]
Now, if in addition $x_1^{n_1}\cdots x_r^{n_r}\equiv 1\pmod {p^l}$ then  $F(t)=G(x_1,\ldots,x_r)\equiv 0\pmod{p^l}$.

It remains to show that $F(t)$ is not the zero polynomial. Notice that if the $a_j$'s are not distinct this is not always true
 (e.g., for $r=2$ if $a_1=a_2$ and $n_1=n_2$ then $F=0$). To show that in our case $F$ does not vanish, we think of it as a complex valued polynomial and note that for it to be identically zero there has to be some choice of signs $\sigma\in\{\pm 1\}^r$ so that the function
\[G_\sigma(t)=\prod_{j=1}^r(\frac{1}{2a_j}(t+\sqrt{t^2-4a_j^2})^{\sigma_j n_j}\]
satisfies $G_\sigma(t)\equiv1$. Assume that there is such a choice $\sigma$, so the derivative $G_\sigma'(t)$ must vanish for all $t$. But on the other hand for such $G_\sigma$ we have
\[G_\sigma'(t)=\sum_{j=1}^r\frac{\sigma_jn_j}{\sqrt{t^2-4a_j^2}},\]
so as $t\to\pm2 a_j$ the term $\frac{\sigma_jn_j}{\sqrt{t^2-4a_j^2}}$ blows up while the rest of the terms remain bounded (here we use that the $a_j$'s are distinct). In particular $G'_\sigma(t)$ is not identically zero.
\end{proof}

% %-------------------------
% %GATHER{/home/member/kelmerdu/My_Documents/Tex/Bib/Mybib}   % For Gather Purpose Only
% \bibliographystyle{amsplain}
% \bibliography{/home/member/kelmerdu/My_Documents/Tex/Bib/Mybib.bib}
% %*********************************************************************
\def\cprime{$'$}
\providecommand{\bysame}{\leavevmode\hbox to3em{\hrulefill}\thinspace}
\providecommand{\MR}{\relax\ifhmode\unskip\space\fi MR }
\providecommand{\MRhref}[2]{%
  \href{http://www.ams.org/mathscinet-getitem?mr=#1}{#2}
}
\providecommand{\href}[2]{#2}

%********************************************************************
\end{document}